\documentclass[11pt,twoside,a4paper]{article}
\usepackage{amsmath,amssymb,amsthm}

\usepackage{graphicx,psfrag}

\newtheorem{Thm}{Theorem}[section]
\newtheorem{Lem}[Thm]{Lemma}
\newtheorem{Pro}[Thm]{Proposition}
\newtheorem{Cor}[Thm]{Corollary}

\theoremstyle{definition}

\newtheorem{Exas}[Thm]{Examples}

\theoremstyle{remark}
\newtheorem{Rem}[Thm]{Remark}

\newcommand{\R}{\mathbb{R}}

\newcommand{\N}{\mathbb{N}}

\newcommand{\cC}{\mathcal{C}}

\newcommand{\cN}{\mathcal{N}}

\newcommand{\cR}{\mathcal{R}}

\newcommand{\cU}{\mathcal{U}}
\newcommand{\cV}{\mathcal{V}}
\newcommand{\cW}{\mathcal{W}}

\newcommand{\al}{\alpha}

\newcommand{\ga}{\gamma}
\newcommand{\Ga}{\Gamma}
\newcommand{\de}{\delta}
\newcommand{\De}{\Delta}
\newcommand{\ep}{\varepsilon}
\newcommand{\om}{\omega}

\newcommand{\si}{\sigma}
\newcommand{\Si}{\Sigma}

\newcommand{\la}{\lambda}
\newcommand{\La}{\Lambda}
\renewcommand{\phi}{\varphi}

\newcommand{\rank}{\operatorname{rank}}
\newcommand{\subcorank}{\operatorname{corank}}
\newcommand{\dist}{\operatorname{dist}}
\newcommand{\diam}{\operatorname{diam}}

\newcommand{\hyp}{\operatorname{H}}

\newcommand{\Lip}{\operatorname{Lip}}

\newcommand{\const}{\operatorname{const}}

\newcommand{\hypdim}{\operatorname{hypdim}}
\newcommand{\asdim}{\operatorname{asdim}}
\newcommand{\sh}{\operatorname{sh}}
\newcommand{\ba}{\operatorname{ba}}
\newcommand{\st}{\operatorname{st}}

\newcommand{\cone}{\operatorname{Co}}
\newcommand{\an}{\operatorname{An}}

\newcommand{\es}{\emptyset}
\renewcommand{\d}{\partial}
\newcommand{\di}{\d_{\infty}}
\newcommand{\set}[2]{\{#1:\,\text{#2}\}}
\newcommand{\sm}{\setminus}
\newcommand{\sub}{\subset}

\newcommand{\ov}{\overline}

\newcommand{\anglin}{\angle_{\infty}}

\begin{document}

\title{Hyperbolic dimension of metric spaces}
\author{Sergei Buyalo\footnote{Supported by RFFI Grant
02-01-00090, CRDF Grant RM1-2381-ST-02 and
SNF Grant 20-668 33.01}
\ \& Viktor Schroeder\footnote{Supported by Swiss National
Science Foundation}}

\date{}
\maketitle

\begin{abstract} We introduce a new quasi-isometry invariant of
metric spaces called the hyperbolic dimension,
$\hypdim$,
which is a version
of the Gromov's asymptotic dimension,
$\asdim$.
One always has
$\hypdim\le\asdim$,
however, unlike the asymptotic dimension,
$\hypdim\R^n=0$
for every Euclidean space
$\R^n$
(while
$\asdim\R^n=n$).
This invariant possesses usual properties of
dimension like monotonicity and product theorems.
Our main result says that
the hyperbolic dimension of any Gromov hyperbolic space
$X$
(with mild restrictions)
is at least the topological dimension of the boundary at infinity
plus 1,
$\hypdim X\ge\dim\di X+1$.
As an application we obtain that there is no quasi-isometric
embedding of the real hyperbolic space
$\hyp^n$
into the metric product of
$n-1$
metric trees stabilized by any Euclidean factor,
$T_1\times\dots\times T_{n-1}\times\R^m$, $m\ge 0$.
\end{abstract}

\section{Introduction}

We introduce a new quasi-isometry invariant of
metric spaces called the hyperbolic dimension,
$\hypdim$,
which is a version
of the Gromov's asymptotic dimension,
$\asdim$.
One always has
$\hypdim\le\asdim$,
however, unlike the asymptotic dimension,
$\hypdim\R^n=0$
for every Euclidean space
$\R^n$
(while
$\asdim\R^n=n$).
This invariant possesses usual properties of
dimension like monotonicity and product theorems. To
formulate our main result, we recall that a metric space
$X$
has {\em bounded growth at some scale}, if for some
constants
$r$, $R$ with
$R>r>0$,
and
$N\in\N$
every ball of radius
$R$
in
$X$
can be covered by
$N$
balls of radius
$r$,
see \cite{BoS}.

\begin{Thm}\label{Thm:main} Let
$X$
be a geodesic Gromov hyperbolic space, which has bounded
growth at some scale and whose boundary at infinity
$\di X$
is infinite. Then
$$\hypdim X\ge\dim\di X+1.$$
\end{Thm}

As an application we obtain.

\begin{Cor}\label{Cor:nonemb} For every
$n\ge 2$
there is no quasi-isometric embedding
$\hyp^n\to T_1\times\dots\times T_{n-1}\times\R^m$
of the real
$n$-dimensional
hyperbolic space
$\hyp^n$
into the product of
$n-1$
trees stabilized by any Euclidean factor
$\R^m$, $m\ge 0$.
\end{Cor}

\begin{Rem} For
$n=2$
this result has been proved in \cite{BS2} by
a different method.
\end{Rem}

\begin{Rem} In \cite{BS2} we have constructed for every
$n\ge 2$
a quasi-isometric embedding of
$\hyp^n$
into the
$n$-fold
product of homogeneous trees whose vertices have an infinite
(countable) degree and whose edges have length 1. By
Corollary~\ref{Cor:nonemb}, that embedding is optimal
with respect to the number of tree-factors even if we
allow the stabilization by Euclidean factors. Furthermore,
it follows that
$\hypdim\hyp^n=n$
for every
$n\ge 2$.

In \cite{JS} it was shown that for every
$n$
there exists a right angled Gromov hyperbolic Coxeter group
$\Ga_n$
with virtual cohomological dimension and coloring number
equal to
$n$.
In \cite{DS}, a bilipschitz embedding
$f_n:X_n\to T\times\dots\times T$
of the Cayley graph
$X_n$
of
$\Ga_n$
into the
$n$-fold
product of an arbitrary exponentially branching tree
$T$
has been constructed. The Cayley graph of any Gromov hyperbolic
group satisfies the conditions of Theorem~\ref{Thm:main}, and
$\dim\di X_n=\dim\di\Ga_n=n-1$
by a result from \cite{BM}(every finitely generated Coxeter group
is virtually torsion free). Thus Theorem~\ref{Thm:main} implies
(see Theorem~\ref{Thm:nonemb} below) that
the embedding
$f_n$
is optimal w.r.t. the number of tree-factors even if we allow
the stabilization by Euclidean factors.
\end{Rem}

\begin{center}{\bf\large Hyperbolic dimension versus subexponential corank}
\end{center}

In our earlier paper \cite{BS1} we introduced another quasi-isometry
invariant of metric spaces called subexponential corank. This invariant
gives an upper bound for the topological dimension of a Gromov hyperbolic
space which can be quasi-isometrically embedded into a given metric space
$X$, $\rank_h(X)\le\subcorank(X)$.
Thus
$\subcorank$
is a useful tool for finding obstacles to such embeddings, and
it works perfectly well in many cases, see \cite{BS1} for details.
However, not in all, e.g., for quasi-isometric embeddings
$\hyp^n\to T_1\times\dots\times T_k\times\R^m$
it gives only
$k\ge n-2$,
while
$\hypdim$
gives optimal
$k\ge n-1$
by Corollary~\ref{Cor:nonemb}. This drawback of
$\subcorank$
is closely related to that
$\subcorank(T)=1>0=\dim\di T$
for every metric tree
$T$,
while
$\subcorank(X)=\dim\di X$
for every CAT($-1$) Hadamard manifold
$X$.

On the other hand,
$\hypdim$,
which is perfect for product of trees, is much harder
to compute than
$\subcorank$.
For example, we do not even know the precise value of
$\hypdim(\hyp^2\times\hyp^2)$
(it must be 3 or 4). We also do not see any direct way to compare
$\subcorank$
and
$\hypdim$.
At the present stage of knowledge, it looks like that these
two invariants are in a sense independent, and each of
them works perfectly well in its own range while failing
in the other.

\medskip\noindent
{\bf Structure of the paper.} In section~\ref{sect:ubg} we introduce
and discuss properties of a class of metric spaces with uniformly
bounded growth rate, UBG-spaces, which is a key ingredient of
the hyperbolic dimension. The main result here is
Proposition~\ref{Pro:endubg}, which is an important step in the proof
of Theorem~\ref{Thm:main}.

In section~\ref{sect:threedef} we give three definitions of the
hyperbolic dimension following the standard line of the
topological dimension theory, and prove their equivalence.
It is convenient to use different definitions in different
situations. In section~\ref{sect:properties} we discuss
properties of the hyperbolic dimension and prove monotonicity
and product theorems. Section~\ref{sect:proof} is devoted
to the proof of Theorem~\ref{Thm:main} and Corollary~\ref{Cor:nonemb}.

Here we briefly recall some notions which are used in the body of
the paper. We denote by
$|x-x'|$
the distance in a metric space
$X$
between
$x$, $x'\in X$.
A map
$f:X\to Y$
between metric spaces is quasi-isometric if for some
$\La\ge 1$, $\la\ge 0$
the estimates
$$\frac{1}{\La}|x-x'|-\la\le|f(x)-f(x')|\le\La|x-x'|+\la,$$
hold for every
$x$, $x'\in X$.
In this case we also say that
$f$
is
$(\La,\la)$-quasi-isometric.
A metric space is geodesic if every two its points
are connected by a geodesic. By a CAT($-1$)-space we
mean a complete, geodesic space whose triangles are
thinner than the comparison triangles in the real
hyperbolic plane
$\hyp^2$.

\medskip\noindent
{\bf Acknowledgment.} The first author is pleased to acknowledge
the hospitality and the support of the University of Z\"urich
where this research has been carried out.

\section{Spaces with uniformly bounded growth rate}\label{sect:ubg}

Here we introduce a class of metric spaces which is a key
ingredient of the notion of the hyperbolic dimension.

\subsection{Definition and properties}
We say that a metric space
$X$
has uniformly bounded growth rate (or is an UBG-space)
if for every
$\rho>0$
there exist
$N\in\N$
and
$r_0>0$
so that every ball of radius
$r\ge r_0$
in
$X$
contains at most
$N$
points which are
$\rho r$-separated.

Equivalently,
$X$
is UBG if for every
$\si>1$
there exist
$N$
and
$r_0$
such that every ball of radius
$\si r$
in
$X$
with
$r\ge r_0$
can be covered by
$N$
balls of radius
$r$.
Or
$X$
is UBG if for every
$\rho>0$
there is
$N$
such that every ball of radius 1 in the scaled
$\la X$
for all sufficiently small
$\la>0$
contains at most
$N$
points which are
$\rho$-separated.

\begin{Rem} The notion of UBG-spaces is close to the notion
of {\em asymptotically doubling} spaces,
where a metric space
$X$
is asymptotically doubling if for some positive
constants
$N$
and
$r_0$
every ball in
$X$
of radius
$r\ge r_0$
can be covered by at most
$N$
balls of radius
$r/2$.
One can show essentially by the same argument as in
\cite[p.\,295]{BoS} that if a geodesic space
$X$
is asymptotically doubling then it is UBG. However,
the assumption that
$X$
is geodesic is too restrictive for our purposes. UBG-spaces
typically appear in our work as preimages of some sets
under quasi-isometric maps, and there is no reason for
them to be geodesic. In other words, we study and use UBG-spaces
as tools rather than in their own right.
\end{Rem}

For technical reason, it is convenient to characterize an
UBG-space by two functions. We let
$\cN$
be the set of functions
$N:(0,1)\to\N$
and
$\cR$
be the set of functions
$R:(0,1)\to[0,\infty)$.
Then the definition above says that
$X$
is
UBG if for some functions
$N\in\cN$
and
$R\in\cR$
and for every
$\rho\in(0,1)$
every ball of radius
$r\ge R(\rho)$
in
$X$
contains at most
$N(\rho)$
points which are
$\rho r$-separated.
In this case we say that
$X$
is
$(N,R)$-bounded.
We also say that
$X$
is
$N$-bounded
if it is
$(N,R)$-bounded
for some
$R\in\cR$.

\begin{Exas} 1. Any Euclidean space
$\R^n$
is
$(N,R)$-bounded
for
$N(\rho)\asymp\rho^{-n}$
and
$R(\rho)\equiv 0$.

The basic example of UBG-spaces is this.

\noindent
2. Let
$B$
be a bounded metric space. Then the metric product
$X=B\times\R^n$
is
$(N,R)$-bounded
for
$N(\rho)\asymp\rho^{-n}$
and
$R(\rho)\ge\frac{2\diam B}{\rho}$.
We emphasize that in this example the function
$N$
counting the number of separated points is actually
independent of
$B$,
while the function
$R$
describing the corresponding scales
tends to infinity as
$\diam B\to\infty$
if one takes as
$B$,
say, an
$\R$-tree.
\end{Exas}

Two functions
$N=N(\rho)$
and
$R=R(\rho)$
are included in the definition of UBG-spaces
instead of fixing some
$\rho\in(0,1)$
for the purpose to make this notion quasi-isometry invariant.

\begin{Lem}\label{Lem:indubg} Let
$f:X\to Y$
be a
$(\La,\la)$-quasi-isometric
map, where
$\La\ge 1$, $\la\ge 0$.
Assume that
$Y$
is
$(N,R)$-bounded
for some
$N\in\cN$
and
$R\in\cR$.
Then
$X$
is
$(N',R')$-bounded
for
$N'(\rho)=N(\rho/2\La^2)$
and
$R'(\rho)=\max\{\frac{1}{\La}(R(\rho)-\la),
  \frac{\la}{\La}(1+\frac{1}{\rho})\}$.
\end{Lem}

\begin{proof} Fix
$\rho\in (0,1/2\La^2)$.
Then for all
$r\ge R(\rho)$
every ball
$B_r\sub Y$
contains at most
$N(\rho)$
$\rho r$-separated
points. We have
$f(B_{r'})\sub B_{\La r'+\la}$
for every ball
$B_{r'}\sub X$,
hence, for
$r'\ge\frac{1}{\La}(R(\rho)-\la)$
the set
$f(B_{r'})$
contains at most
$N(\rho)$
$\si$-separated
points,
$\si=\rho(\La r'+\la)$.
Take
$\si'=\La(\si+\la)$.
If
$x$, $x'\in X$
are
$\si'$-separated,
then
$f(x)$, $f(x')\in Y$
are
$\si$-separated.
It follows that the ball
$B_{r'}$
itself contains at most
$N(\rho)$
$\si'$-separated
points.
We put
$R'(\rho)=\max\{\frac{1}{\La}(R(\rho)-\la),
  \frac{\la}{\La}(1+\frac{1}{\rho})\}$,
$\rho'=2\La^2\rho$.
Hence, for
$r'\ge R'(\rho)$
we have
$\si'\le 2\La^2\rho r'=\rho'r'$,
and
$\rho'r'$-separated
points are certainly
$\si'$-separated.
Then for every
$r'\ge R'(\rho)$
every ball
$B_{r'}\sub X$
contains at most
$N(\rho)$
$\rho'r'$-separated
points, and the space
$X$
is
$(N',R')$-bounded,
where
$N'(\rho)=N(\rho/2\La^2)$
and
$R'$
as above.
\end{proof}

\begin{Cor}\label{Cor:ubgqi} The property of a metric
space to have a uniformly bounded growth rate is a
quasi-isometry invariant.
\qed
\end{Cor}

\begin{Lem}\label{Lem:ubgpol} Every UBG space
$Y$
has a polynomial growth rate, that is, there exists
$k=k(Y)>0$
such that for every (sufficiently large)
$\de$
every ball of radius
$r$
in
$Y$
contains at most
$d\cdot r^k$
$\de$-separated
points provided
$r$
is sufficiently large, where the constant
$d>0$
depends only on
$Y$
and
$\de$.
\end{Lem}

\begin{proof} We fix
$\rho\in(1/2,1)$
and take
$\de$
so that
$\de'=\frac{2\rho-1}{\rho}\de\ge R(\rho)$.
Then
every ball
$B_{\de'}\sub Y$
of radius
$\de'$
contains at most
$N=N(\rho)$
points which are
$\de'$-separated.

Take
$r\ge\de$.
Then
$r>R(\rho)$,
and any ball
$B_r\sub Y$
contains at most
$N$
points which are
$\rho r$-separated.
Take a maximal
$\rho r$-separated
subset in
$B_r$.
Then the balls
$B_{\rho r}$
centered at its points cover
$B_r$,
and their number is at most
$N$.
Applying this argument to every
$B_{\rho r}$
from the covering of a fixed
$B_r$
and proceeding by induction, we obtain that
$B_r$
contains at most
$N^q$
points which are
$\rho^qr$-separated,
provided
$\rho^qr\ge\de'$.
There is
$q\in\N$
with
$\rho^qr<\de\le\rho^{q-1}r$.
For this
$q$
the condition
$\rho^qr\ge\de'$
is fulfilled, thus every
$B_r$
contains at most
$N^q\le d\cdot r^k$
points which are
$\de$-separated, where
$k=\ln N/\ln\frac{1}{\rho}$,
$d=(\rho\de)^{-k}$.
\end{proof}

\begin{Lem}\label{Lem:ubgprod} Let
$X$, $Y$
be UBG-spaces. Then
$Z=X\times Y$
is an UBG-space.
\end{Lem}

\begin{proof} The proof is straightforward, if one uses
the covering definition of UBG.
\end{proof}

\subsection{UBG-subsets in a CAT($-1$)-space}\label{subsect:catubg}

Let
$X$
be a CAT($-1$)-space. We fix a base point
$x_0\in X$
and define the angle metric
$\anglin$
in the boundary at infinity
$\di X$
as follows. Given
$\xi$, $\xi'\in\di X$,
we consider the unit speed geodesic rays
$c_{\xi}$, $c_{\xi'}$
from
$x_0$
to
$\xi$, $\xi'$
respectively, and put
$$\anglin(\xi,\xi')=\lim_{s\to\infty}
  \angle(\ov{c}_{\xi}(s)\ov o\ \ov{c}_{\xi'}(s)),$$
where
$\angle(\ov{c}_{\xi}(s)\ov o\ \ov{c}_{\xi'}(s))$
is the angle at
$\ov o$
of the comparison triangle in
$\hyp^2$
for the triangle
$c_{\xi}(s)x_0c_{\xi'}(s)$.
From the hyperbolic geometry we know that
$$\tan\left(\frac{1}{4}\anglin(\xi,\xi')\right)=
  e^{-\dist(\ov o,\ov\xi\,\ov\xi')},$$
where
$\ov\xi$, $\ov\xi'\in\di\hyp^2$
satisfy
$\angle_{\ov o}(\ov\xi,\ov\xi')=\anglin(\xi,\xi')$,
and
$\ov\xi\,\ov\xi'$
is the geodesic in
$\hyp^2$
with the end points at infinity
$\ov\xi$, $\ov\xi'$.
Thus
$\anglin(\xi,\xi')\le 4e^{-\dist(\ov o,\ov\xi\,\ov\xi')}$.

{\em The shadow} of a set
$A\sub X$
is a subset
$\sh(A)\sub\di X$
which consists of the ends
$\xi$
of all rays
$x_0\xi$
intersecting
$A$
(so
$\sh(x_0)=\di X$).
Given
$\de>0$
we define {\em the angle
$\de$-measure}
of
$A$,
$\angle_\de A$,
as
$$\angle_\de A=\inf_{\cC}\sum_{B\in\cC}\diam(\sh(B)),$$
where the infimum is taken over all coverings
$\cC$
of
$A$
by balls of radius
$\ge\de$
in
$X$.

\begin{Lem}\label{Lem:ubgshadow} Given functions
$N\in\cN$, $R\in\cR$,
for every sufficiently large
$\de$
there is a positive constant
$C$
depending only on
$\rho\in(1/2,1)$, $N(\rho)$, $R(\rho)$
and
$\de$
such that
if a subset
$A\sub X$
is
$(N,R)$-bounded
and
$\dist(x_0,A)\ge c>\de$,
then
$$\angle_\de A\le C\cdot e^{-c/2},$$
\end{Lem}

\begin{proof} As in the proof of Lemma~\ref{Lem:ubgpol}, we fix
$\rho\in(1/2,1)$
and take
$\de\ge\frac{\rho}{2\rho-1}R(\rho)$. Then, by Lemma~\ref{Lem:ubgpol},
every ball
$B_r\sub X$
with
$r>\de$
contains at most
$dr^k$
points of
$A$
which are
$\de$
separated, where
$k$
depends on
$\rho$, $N(\rho)$,
and
$d=(\rho\de)^{-k}$.
Furthermore, since
$k$
is independent of
$\de$,
we can assume that
$e^{c/2}\ge(c+1)^k$
for each
$c>\de$.

Take a maximal
$\de$-separated
subset
$A'\sub A$.
Then
$A\sub\cup_{a\in A'}B_\de(a)$.
For any ball
$B_\de(a)$, $a\in A'$,
consider
$\xi$, $\xi'\in\sh(B_\de(a))$
with
$\anglin(\xi,\xi')=\diam(\sh(B_\de(a)))$.
Then
$\diam(\sh(B_\de(a)))\le 4e^{-\dist(\ov o,\ov\xi\,\ov\xi')}$
in notation introduced above. We take
$x\in x_0\xi\cap B_\de(a)$, $x'\in x_0\xi'\cap B_\de(a)$
and consider the piecewise geodesic curve
$\ga$
in
$X$,
which consists of the geodesic rays
$x\xi$, $x'\xi'$
and the segment
$xx'$.
The curve
$\ga$,
as well as the geodesic
$\xi\xi'$,
connects in
$X$
the points
$\xi$, $\xi'$,
and
$\dist(x_0,\ga)\ge\dist(x_0,B_\de(a))=|x_0-a|-\de$.
Furthermore,
$\dist(x_0,\ga)\le\dist(\ov o,\ov\xi\,\ov\xi')$
by comparison with
$\hyp^2$.
Thus
$$\diam(\sh(B_\de(a)))\le 4e^{\de-|x_0-a|}$$
and
$$\angle_\de A\le\sum_{a\in A'}\diam(\sh(B_\de(a))).$$
Since
$c>\de$,
for each
$\tau\ge c+1$,
the number of points from
$A'$
whose distances to
$x_0$
lie in the interval
$[\tau-1,\tau)$
is
$\le d\cdot\tau^k$.
Thus we have
\begin{eqnarray*}
 \angle_\de A&\le& 4e^{\de}\sum_{a\in A'}e^{-|x_0-a|}\le
   4de^{\de}\left(\sum_{q=0}^{\infty}(c+q+1)^ke^{-q}\right)e^{-c}\\
   &\le&4de^{\de}\left(\sum_{q=0}^{\infty}(q+1)^ke^{-q}\right)(c+1)^ke^{-c}
   \le C\cdot e^{-c/2}
\end{eqnarray*}
by the choice of
$c$.
\end{proof}

\subsection{A cut property of UBG-subsets}

A subset
$A$
of a metric space
$X$
is said to be roughly connected if for some
$\si>0$
and for every
$a$, $a'\in A$
there is a sequence
$a_0=a,\dots,a_k=a'$
in
$A$
with
$|a_i-a_{i-1}|\le\si$, $i=1,\dots,k$.
Such a sequence is called a rough or a
$\si$-path
between
$a$
and
$a'$,
and we also say that
$A$
is
$\si$-connected.
One says that a roughly connected subset
$A$
of a geodesic space
$X$
is {\em cut-quasi-convex,}
if there is
$c>0$
such that for every
$a$, $a'\in A$
and every
$x\in aa'$,
every rough path in
$A$
between
$a$, $a'$
intersects the ball
$B_c(x)\sub X$
of radius
$c$
centered at
$x$.
In this case we also say that
$A$
is
$c$-cut-convex,
and
$c$
is called {\em the cut radius} of
$A$.
This property, obviously, implies that every geodesic
segment
$aa'\sub X$
with the end points
$a$, $a'\in A$
lies in the
$c$-neighborhood
of
$A$,
i.e., the set
$A$
is in particular quasi-convex. This justifies our terminology.

\begin{Pro}\label{Pro:cutubg} Assume that a
$\si$-connected
subset
$A$
of CAT($-1$)-space
$X$
is
$(N,R)$-bounded
for some functions
$N\in\cN$, $R\in\cR$.
Then
$A$
is cut-quasi-convex, and the cut radius
$c$
depends only on
$\rho\in(1/2,1)$, $N(\rho)$, $R(\rho)$
and
$\si$,
$c=c(\rho,N(\rho),R(\rho),\si)$.
\end{Pro}

\begin{proof} Fix
$\rho\in(1/2,1)$
and take a sufficiently large
$\de$
provided by Lemma~\ref{Lem:ubgshadow}. Furthermore,
we can assume that
$\de\ge\si$.
Next, we take
$c'>\de$
such that
$C\cdot e^{-c'/2}<\pi$,
where
$C$
is the constant from Lemma~\ref{Lem:ubgshadow},
and put
$c=c'+\de$.
Then
$c=c(\rho,N(\rho),R(\rho),\si)$.

Assume that for some
$a$, $a'\in A$
there is a
$\si$-path
$\ga$
in
$A$
between
$a$, $a'$
which misses the ball
$B_c(x_0)\sub X$
for some
$x_0\in aa'$.
Then
$A'=\cup_{b\in\ga}B_\de(b)$
is a connected subset in
$X$
containing the points
$a$, $a'$,
thus
$\angle_\de A'\ge\pi$
for every
$\de>0$.

On the other hand,
$\dist(x_0,A')\ge c'$,
and we have
$\angle_\de A\le C\cdot e^{-c'/2}<\pi$
by Lemma~\ref{Lem:ubgshadow}. This is a contradiction,
and hence
$A$
is cut-quasi-convex with the cut radius
$c$.
\end{proof}

\begin{Cor}\label{Cor:hypubg} If a ball
$B_r$
of radius
$r$
in a CAT($-1$)-space
$X$
is
$(N,R)$-bounded
for some
$N\in\cN$, $R\in\cR$,
then
$r\le c$
for some constant
$c=c(\rho,N(\rho),R(\rho))$, $\rho\in(1/2,1)$.
\qed
\end{Cor}

Discussing the cut-quasi-convex property, we have used so far only
the fact that any UBG-space has a polynomial growth rate (actually,
a subexponential growth rate suffices). In what follows, the next
Proposition plays a key role, and we use in it the whole power of
the definition of UBG-spaces.

\begin{Pro}\label{Pro:endubg} Let
$A$
be an
$N$-bounded
(for some function
$N\in\cN$),
subset in a CAT($-1$)-space
$X$.
Then for every
$\si>0$
the union
$\di A_\si$
of the boundaries at infinity of
$\si$-connected
components of
$A$
contains at most
$M$
points, where
$M<\infty$
depends only on the function
$N$.
\end{Pro}

\begin{proof} Fix
$x_0\in X$.
It follows from the cut-quasi-convex property of
any roughly connected component of
$A$
that a tail of every geodesic ray
$x_0\xi\sub X$, $\xi\in\di A_\si$,
lies in the
$c$-neighborhood
of
$A$
for some
$c>0$
depending only on the bounding parameters for
$A$
and
$\si$.
Thus for every
$\rho\in(0,1)$
and for all sufficiently large
$r$
the ball
$B_r(x_0)\sub X$
contains at least as much
$\rho r$-separated
points from
$A$
as the cardinality of
$\di A_\si$.
Hence, the claim.
\end{proof}

\section{Three definitions of the hyperbolic dimension}\label{sect:threedef}

The hyperbolic dimension is a variation of the Gromov's asymptotic
dimension (see \cite[1.E]{Gr}) with only difference that we take as
``small'' the UBG-sets. Here we give three equivalent definitions
of the hyperbolic dimension following the standard line of topological
dimension theory. As in \cite{BD}, we find convenient to
use the Lebesgue number of a covering in our definitions instead of
$d$-multiplicity
as in \cite[1.E]{Gr}.

Recall that the Lebesgue number of a covering
$\cU$
of a metric space
$X$
is the maximal radius
$L(\cU)$
such that any (open) ball in
$X$
of that radius is contained in some element of the covering,
$$L(\cU)=\inf_{x\in X}\max\set{\dist(x,X\sm U)}{$U\in\cU$}.$$

Given
$N\in\cN$,
a covering
$\cU$
of a metric space
$X$
is said to be {\em uniformly
$N$-bounded,}
if
\begin{itemize}
\item{} there is a function
$R\in\cR$
such that every element of the covering is
$(N,R)$-bounded;
\item{} any finite union of elements of the covering is
$N$-bounded.
\end{itemize}

\begin{Rem}\label{Rem:boundcover} To get a better grip on
the second property, which is rather strong,
assume that
$X$
is a CAT($-1$)-space. Then by Proposition~\ref{Pro:endubg}, the
boundary at infinity
$\di A_{\si}\sub\di X$
of roughly connected components of any finite union
$A=\cup_iU_i$
of elements
$U_i\in\cU$
has the cardinality bounded above independently of
the number of the elements.
\end{Rem}

\subsection{First definition}

The hyperbolic dimension of a metric space
$X$, $\hypdim_1 X$,
is the minimal
$n$
such that for every
$d>0$
there are a function
$N\in\cN$
and a covering of
$X$
by
$n+1$
subsets
$X_j=\cup_{\al}X_{j\al}$, $j=1,\dots,n+1$
such that
\begin{itemize}
\item{} $X_{j\al}\cap X_{j\al'}=\es$
for every
$j=1,\dots,n+1$
and all
$\al\neq\al'$;
\item{} the covering
$\{X_{j\al}\}$
of
$X$
is uniformly
$N$-bounded
and its Lebesgue number is
$\ge d$.
\end{itemize}

\subsection{Second definition}
The hyperbolic dimension of a metric space
$X$, $\hypdim_2 X$,
is the minimal
$n$
such that for every
$d>0$
there are a function
$N\in\cN$
and a covering
$\cU$
of
$X$
with Lebesgue number
$L(\cU)\ge d$
and multiplicity
$\le n+1$,
which is uniformly
$N$-bounded.

Clearly,
$\hypdim_2X\le\hypdim_1X$.

\subsection{Third definition}

Given a index set
$J$,
we let
$R^J$
be the Euclidean space of functions
$J\to\R$
with finite support, i.e.,
$x\in\R^J$
iff only finitely many coordinates
$x_j=x(j)$
are not zero. Then the distance is well defined by
$$|x-x'|^2=\sum_{j\in J}(x_j-x_j')^2.$$
Let
$\De^J\sub\R^J$
be the standard simplex, i.e.,
$x\in\De^J$
iff
$x_j\ge 0$
for all
$j\in J$
and
$\sum_{j\in J}x_j=1$.

A metric in
$n$-dimensional
simplicial complex
$P$
is said to be {\em uniform} if
$P$
is isometric to a subcomplex of
$\De^J\sub\R^J$
for some index set
$J$.
Every simplex
$\si\sub P$
is then isometric to the standard
$k$-simplex
$\De^k\sub\R^{k+1}$, $k=\dim\si$
(so, for a finite
$J$,
$\dim\De^J=|J|-1$).
For every simplicial polyhedron
$P$
there is the canonical embedding
$u:P\to\De^J$,
where
$J$
is the vertex set of
$P$,
which is affine on every simplex. Its image
$P'=u(P)$
is called {\em the uniformization} of
$P$,
and
$u$
is the uniformization map.

The hyperbolic dimension of a metric space
$X$, $\hypdim_3X$,
is the minimal
$n$
such that for every
$\la>0$
there are a function
$N\in\cN$
and a
$\la$-Lipschitz
map
$p:X\to P$
into a uniform
$n$-dimensional simplicial
polyhedron
$P$,
for which the covering
$\set{p^{-1}(\st_v)}{$v\in P$}$
of
$X$
by preimages of the open stars
$\st_v\sub P$
of the vertices of
$P$
is uniformly
$N$-bounded.

\subsection{Equivalence of the definitions}

\begin{Pro}\label{Pro:equivhypdim} For every metric space
$X$
we have
$$\hypdim_1X=\hypdim_2X=\hypdim_3X.$$
\end{Pro}

\begin{proof} We have already mentioned that
$\hypdim_2X\le\hypdim_1X$
easily follows from the definitions. The proof of
$\hypdim_3X\le\hypdim_2X$
is fairly standard (see \cite[1.E$_1$]{Gr}, \cite[Propositions~1,2]{BD}).
Denote
$n=\hypdim_2X$.
Given
$d>0$,
we have a function
$N\in\cN$
and a uniformly
$N$-bounded
covering
$\cU=\{U_j\}_{j\in J}$
of
$X$
with multiplicity
$\le n+1$
and Lebesgue number
$L(\cU)\ge d$.
Using these data we construct a
$\la$-Lipschitz
map
$p:X\to\De^J$
with
$\la\le\frac{(n+2)^2}{d}$,
whose image lies in a
$n$-dimensional
subpolyhedron
$P\sub\De^J$,
as follows.

Given
$j\in J$,
we define
$q_j:X\to\R$
by
$q_j(x)=\min\{d,\dist(x,X\sm U_j)\}$.
Then
$\sum_{j\in J}q_j(x)\ge d$
for every
$x\in X$.
Furthermore,
$q_j(x')\le q_j(x)+|x-x'|$
and
$$\sum_{j\in J}q_j(x')\le\sum_{j\in J}q_j(x)+(n+1)|x-x'|,$$
because in each sum there are at most
$n+1$
nonzero summands. Using this one obtains
$$\frac{1}{\sum_{j\in J}q_j(x')}\le\frac{1}{\sum_{j\in J}q_j(x)}
+\frac{(n+1)|x-x'|}{d\cdot\sum_{j\in J}q_j(x)}.$$
Now, we put
$p_j(x)=q_j(x)/\sum_{j\in J}q_j(x)$.
Then abbreviating
$\Si=\sum_{j\in J}q_j(x)$, $\Si'=\sum_{j\in J}q_j(x')$
we obtain
$$p_j(x')-p_j(x)=\frac{q_j(x')-q_j(x)}{\Si'}
  +\left(\frac{1}{\Si'}
  -\frac{1}{\Si}\right)q_j(x)\le\frac{n+2}{d}|x-x'|.$$
Finally, for
$p=\{p_j\}_{j\in J}$
we have
$$|p(x')-p(x)|^2=\sum_{j\in J}(p_j(x')-p_j(x))^2
  \le\frac{(n+2)^2(2n+2)}{d^2}|x-x'|^2,$$
hence
$p:X\to\De^J$
is
$\la$-Lipschitz
with
$\la\le\frac{(n+2)^2}{d}$.
Since for every
$x\in X$
there are at most
$n+1$
nonzero coordinates
$p_j(x)$, $j\in J$,
the image
$p(X)$
lies in a
$n$-dimensional
subpolyhedron
$P\sub\De^J$.
For each vertex
$v\in P$
the preimage of its open star
$p^{-1}(\st_v)\sub X$
is contained in some element of
the covering
$\cU$.
Thus the covering
$\set{p^{-1}(\st_v)}{$v\in P$}$
of
$X$
is uniformly
$N$-bounded.
It follows that
$\hypdim_3X\le\hypdim_2X$.

To prove that
$\hypdim_1X\le\hypdim_3X$,
we assume that for every
$\la>0$
there are a function
$N\in\cN$
and a
$\la$-Lipschitz
map
$p:X\to P$
into a uniform
$n$-dimensional simplicial
polyhedron
$P$,
for which the covering
$\set{p^{-1}(\st_v)}{$v\in P$}$
of
$X$
is uniformly
$N$-bounded.

Every
$j$-dimensional
simplex
$\si\sub P$
is matched by its barycenter, which is the vertex
$v_\si$
in the first barycentric subdivision
$\ba P$
of
$P$.
Let
$P_j$
be the union of the open starts
$P_{j\si}$
of
$\ba P$
of all
$v_\si$
with
$\dim\si=j$,
$$P_j=\bigcup_{\dim\si=j}P_{j\si}.$$
Then
$P_{j\si}\cap P_{j\si'}=\es$
for
$\si\neq\si'$
and
$P=\cup_{j=0}^nP_j$.
Now, we put
$X_{j\si}=p^{-1}(P_{j\si})$, $X_j=p^{-1}(P_j)$.
Then
$X=\cup_{j=0}^nX_j$
and
$X_{j\si}\cap X_{j\si'}=\es$
for every
$j=0,\dots,n$
and all
$\si\neq\si'$.
Furthermore, the stars
$P_{j\si}$
of
$\ba P$
are contained in appropriate open stars of
$P$,
thus the covering
$\{X_{j\si}\}$
of
$X$
is uniformly
$N$-bounded.
Since the polyhedron
$P$
is uniform, there is a lower bound
$l_n>0$
for the Lebesgue number of the covering
$\{P_{j\si}\}$.
Therefore, the Lebesgue number of
$\{X_{j\si}\}$
is at least
$l_n/\la$.
This shows that
$\hypdim_1X\le\hypdim_3X$.
\end{proof}

From now on, we use notation
$\hypdim X$
for the common value
$$\hypdim_1X=\hypdim_2X=\hypdim_3X.$$

\section{Properties of the hyperbolic dimension}\label{sect:properties}

The following two properties of the hyperbolic dimension
are obvious from the definitions:
\begin{itemize}
\item{} $\hypdim X=0$
for every UBG-space
$X$.
\item{} $\hypdim X\le\asdim X$
for every metric space
$X$.
\end{itemize}
(The asymptotic dimension of a metric space
$X$, $\asdim X$,
is defined precisely in the same way as the hyperbolic
dimension taking instead of uniformly
$N$-bound\-ed
coverings, coverings by sets with uniformly bounded diameter,
see \cite{Gr}, \cite{BD}). The last inequality can be strong.
For example,
$\asdim\R^n=n$
for every
$n\ge 0$,
while
$\hypdim\R^n=0$
since
$\R^n$
is an UBG-space.

\subsection{Monotonicity of the hyperbolic dimension}

The following simple but important monotonicity theorem implies, in
particular, that the hyperbolic dimension is a quasi-isometry
invariant of a metric space.

\begin{Thm}\label{Thm:monhypdim}
Let
$f:X\to X'$
be a quasi-isometric map between metric spaces
$X$, $X'$.
Then
$$\hypdim X\le\hypdim X'.$$
\end{Thm}

\begin{proof} The proof is straightforward using
the second definition of the hyperbolic dimension
and Lemma~\ref{Lem:indubg}.
\end{proof}

\begin{Cor}\label{Cor:qinv} The hyperbolic dimension is
a quasi-isometry invariant of metric spaces.
\qed
\end{Cor}

\subsection{Product Theorem for the hyperbolic dimension}
\label{subsect:product}

While simplicial complexes appear as nerves of coverings, they
are not convenient for the proof of the Product Theorem for the
hyperbolic dimension. Instead, we use cubical complexes. Thus
we first describe relation between simplicial and cubical complexes.

Given a index set
$J$,
one defines
$Q^J=\set{x\in\R^J}{$\max_{j\in J}x_j=1,x_j\ge 0$}$.
There is a canonical homeomorphism
$\pi_J:Q^J\to\De^J$
(the radial projection) given by
$$\pi_J(x)=\frac{x}{\sum_{j\in J}x_j}$$
for every
$x\in Q^J$.

\begin{Lem}\label{Lem:lipi} Restricted to any
$n$-dimensional
coordinate subspace
$\R^n\sub\R^J$
the map
$\pi_J$
is Lipschitz with the Lipschitz constant
$\le n+1$.
\end{Lem}

\begin{proof} For
$x$, $x'\in Q^J\cap\R^n$
we have
$|x_j-x_j'|\le|x-x'|$
for
every
$j\in J$
and
$\sum_{j\in J}|x_j-x_j'|\le n|x-x'|$,
$\sum_{j\in J}x_j\ge 1$.
Using this, one easily obtains
$|\pi_J(x)-\pi_J(x')|\le(n+1)|x-x'|$,
and the claim follows.
\end{proof}

The inverse homeomorphism
$\om_J=\pi_J^{-1}:\De^J\to Q^J$
is given by
$\om_J(x)=\la(x)x$,
where
$\la(x)=(\max_{j\in J}x_j)^{-1}$.

\begin{Lem}\label{Lem:lipinv} Restricted to any
$n$-dimensional
coordinate subspace
$\R^n\sub\R^J$
the map
$\om_J$
is Lipschitz with the Lipschitz constant
$\le n(1+\sqrt n)$.
\end{Lem}

\begin{proof} For
$x$, $x'\in\De^J\cap\R^n$
we have
$|x_j-x_j'|\le|x-x'|$
for
every
$j\in J$
and
$$|\max_{j\in J}x_j-\max_{j\in J}x_j'|\le
  \max_{j\in J}|x_j-x_j'|\le|x-x'|,$$
$\max_{j\in J}x_j$, $\max_{j\in J}x_j'\ge 1/n$.
Using this, we easily obtain
$|\om_J(x)-\om_J(x')|\le(n+n\sqrt n)|x-x'|$,
and the claim follows.
\end{proof}

Faces of dimension
$k\ge 0$
of
$Q^J$
are defined as the closures of its subsets where exactly
$k$
coordinates have values in
$(0,1)$.
They are isometric to the unit cube in
$\R^k$,
and every such a face is mapped by
$\pi_J$
onto the union of simplices of the first barycentric
subdivision
$\ba\De^J$.
Vice versa, for every vertex of a
$k$-dimensional
simplex
$\De^k\sub\De^J$
its (closed) star in the first barycentric subdivision
$\ba\De^k$
is mapped by
$\om_J$
homeomorphically onto a
$k$-dimensional
face of
$Q^J$.
Therefore, the image
$\om_J(\De^k)\sub Q^J$
consists of
$k+1$
cubical
$k$-faces.
As a corollary of this and Lemmas~\ref{Lem:lipi}, \ref{Lem:lipinv},
we obtain.

\begin{Lem}\label{Lem:simplcub} For any
$n$-dimensional
subcomplex
$P\sub\De^J$,
the image
$P'=\om_J(P)$
is a
$n$-dimensional
cubical subcomplex of
$Q^J$,
and
$\om_J:P\to P'$
is a Lipschitz homeomorphism with the Lipschitz
constant depending only on
$n$.

Conversely, for any
$n$-dimensional
subcomplex
$P'\sub Q^J$,
the image
$P=\pi_J(P')\sub\De^J$
is a
$n$-dimensional
simplicial subcomplex of
$\ba\De^J$,
and
$\pi_J:P'\to P$
is a Lipschitz homeomorphism with the Lipschitz
constant depending only on
$n$.
\qed
\end{Lem}

Now, we prove the product Theorem for the hyperbolic
dimension.

\begin{Thm}\label{Thm:prodhypdim} For any metric spaces
$X_1$, $X_2$,
we have
$$\hypdim(X_1\times X_2)\le\hypdim X_1+\hypdim X_2.$$
\end{Thm}

\begin{proof} We use the third definition of the hyperbolic
dimension. Given
$\la_k>0$
there are a function
$N_k\in\cN$
and a
$\la_k$-Lipschitz
map
$p_k:X_k\to P_k\sub\De^{J_k}$
into
$n_k$-dimensional
uniform simplicial polyhedron
$P_k$,
where
$n_k=\hypdim X_k$,
such that the covering
$\set{p_k^{-1}(\st_v)}{$v\in P_k$}$
of
$X_k$
is uniformly
$N_k$-bounded,
$k=1,2$.
We assume that
$n_1$, $n_2<\infty$,
since otherwise there is nothing to prove.

By Lemma~\ref{Lem:simplcub},
$P_k'=\om_k(P_k)\sub Q^{J_k}$
is a
$n_k$-dimensional
cubical subcomplex of
$Q^{J_k}$,
and the homeomorphism
$\om_k=\om_{J_k}:P_k\to P_k'$
is Lipschitz with the Lipschitz constant
depending only on
$n_k$, $k=1,2$.
Then
$\om_k\circ p_k:X_k\to P_k'$
is Lipschitz with Lipschitz constant
$\le\const(n_k)\cdot\la_k$.
Furthermore, the covering
$\set{(\om_k\circ p_k)^{-1}(\st_w)}{$w\in P_k'$}$
of
$X_k$
by preimages of the open cubical stars
$\st_w$
of the vertices of
$P_k'$
is uniformly
$N_k$-bounded,
because every such a star
$\st_w\sub P_k'$
lies in
$\om_k(\st_v)$
for an appropriate vertex
$v\in P_k$, $k=1,2$.

Let
$J=J_1\cup J_2$
be the disjoint union. We define
$p:X_1\times X_2\to\R^J$
by
$p(x_1,x_2)=(\om_1\circ p_1(x_1),\om_2\circ p_2(x_2))$
for every
$(x_1,x_2)\in X_1\times X_2$.
Then
$p(X_1\times X_2)\sub P'$,
where
$P'=P_1'\times P_2'\sub Q^J$
is a cubical subcomplex of dimension
$n_1+n_2$,
and the map
$p$
is Lipschitz with Lipschitz constant
$\Lip(p)\le\const(n_1,n_2)\cdot\max\{\la_1,\la_2\}$.
Using Lemma~\ref{Lem:ubgprod} one easily checks that
the covering
$\set{p^{-1}(\st_w)}{$w\in P'$}$
of
$X_1\times X_2$
by preimages of the open cubical stars is uniformly
$N$-bounded
for some function
$N\in\cN$
depending only on
$N_1$, $N_2$.
Applying the homeomorphism
$\pi_J:Q^J\to\De^J$,
we obtain a simplicial subcomplex
$P=\pi_J(P')\sub\ba\De^J$
and a
$\mu$-Lipschitz
map
$\pi_J\circ p:X_1\times X_2\to P$,
$\mu\le\const(n_1,n_2)\cdot\max\{\la_1,\la_2\}$,
with required
$N$-boundedness
property. It remains to compose this map with
a Lipschitz simplicial homeomorphism sending
$P$
into a simplicial subcomplex of
$\De^{J'}$,
where
$J'$
is the set of all nonempty finite subsets in
$J$.
\end{proof}

\section{Proof of Theorem~\ref{Thm:main}}\label{sect:proof}

\subsection{d-multiplicity of a covering}\label{subsect:dmult}

Fix
$d>0$.
Recall that the
$d$-multiplicity
of a covering
$\cU$
of a metric space
$X$
is
$\le n$,
if no ball of radius
$d$
in
$X$
meets more than
$n$
elements of the covering (see \cite[1.E]{Gr}).
One can define an auxiliary hyperbolic dimension
$\hypdim'X$
as a minimal
$n$
such that for every
$d>0$
there are a function
$N\in\cN$
and a uniformly
$N$-bounded
covering
$\cU$
of
$X$
with
$d$-multiplicity
$\le n+1$.

\begin{Lem}\label{Lem:ahypdim} For every metric space
$X$
we have
$\hypdim'X\le\hypdim X$.
\end{Lem}

\begin{proof} Let
$\cU=\{U_j\}_{j\in J}$
be a uniformly
$N$-bounded
covering of
$X$
with multiplicity
$\le n+1$, $n=\hypdim X$,
and Lebesgue number
$L(\cU)\ge 2d$
for some function
$N\in\cN$
and some
$d>0$.
Following \cite[Assertion~1]{BD}, we define
$V_j=U_j\sm D_d(X\sm U_j)$, $j\in J$,
where
$D_d(A)$
is the metric
$d$-neighborhood
of
$A\sub X$.
Since
$L(\cU)\ge 2d$,
the collection
$\cV=\{V_j\}_{j\in J}$
is still a covering of
$X$.
Moreover,
$\cV$
is uniformly
$N$-bounded
because
$V_j\sub U_j$
for every
$j\in J$.
Furthermore, if a ball
$B_d(x)\sub X$
meets some
$V_j$,
then
$x\in U_j$.
Thus the
$d$-multiplicity of
$\cV$
is
$\le n+1$.
Hence, the claim.
\end{proof}

\subsection{Ends of uniformly
$N$-bounded
coverings of
$X$}

Recall that a metric space
$X$
has bounded geometry if there are
$\rho_X>0$
and a function
$M_X:(0,\infty)\to(0,\infty)$
such that every ball
$B_r\sub X$
of radius
$r>0$
contains at most
$M_X(r)$
points which are
$\rho_X$-separated.

\begin{Lem}\label{Lem:boundend} Let
$X$
be a CAT($-1$)-space
with bounded geometry, and let
$\cU$
be a uniformly
$N$-bounded
covering of
$X$
with finite
$d$-multi\-plicity
for a function
$N\in\cN$
and a sufficiently large
$d$.
Furthermore, assume that the elements of
$\cU$
are
$\si$-connected
for some
$\si\ge 10d$.
Then
$\di U\sub\di X$
is finite for every
$U\in\cU$
and the number of different elements of
$\cU$
with infinite diameter is finite.
\end{Lem}

\begin{proof} Since the covering
$\cU$
is uniformly
$N$-bounded,
there is a function
$R\in\cR$
such that every element of the covering is
$(N,R)$-bounded.
Then by Proposition~\ref{Pro:cutubg} every
$U\in\cU$
is cut-quasi-convex with cut radius
$c>0$
depending only on
$N$, $R$
and
$\si$.
By Proposition~\ref{Pro:endubg} there is
an upper bound
$M<\infty$
for the cardinality of the union
$\di A_\si$
of the boundaries
at infinity of
$\si$-connected
components
of any
$N$-bounded
$A\sub X$.
In particular,
$|\di U|\le M$
for every
$U\in\cU$.

We assume that
$X$
has
$(\rho_X,M_X)$-bounded
geometry, and that
$d$-multi\-plicity
of
$\cU$
is
$\le n$
for
$d\ge\rho=\rho_X$.
Then we put
$M_0=(n+1)\cdot M\cdot M_X(2c)$
and assume that there are
$\ge M_0$
different elements of the covering
$\cU$
with infinite diameter. Since every finite
union of elements of
$\cU$
is
$N$-bounded,
there is
$\xi\in\di X$
which is a common point of
$\di U$
for at least
$M_0/M$
different elements
$U\in\cU$.

We fix
$x_0\in X$
and consider the geodesic ray
$x_0\xi\sub X$.
By the cut-quasi-convex property of
$U$,
a tail
$x_1\xi\sub x_0\xi$
lies in the
$c$-neighborhood
of every
selected above
$U$.
Thus for every
$x\in x_1\xi$
the ball
$B_{2c}(x)$
intersects at least
$(n+1)\cdot M_X(2c)$
different elements of the covering. On the other hand,
$B_{2c}(x)$
can be covered by
$\le M_X(2c)$
balls of radius
$\rho$.
Hence, there is a ball
$B_\rho(x')\sub X$
which intersects at least
$n+1$
different elements of the covering. This is
a contradiction, because
$d\ge\rho$
and
$d$-multiplicity
of the covering is
$\le n$.
Thus there is at most
$M_0-1$
different elements of
$\cU$
with infinite diameter.
\end{proof}

\subsection{Radial contraction}

\begin{Lem}\label{Lem:coveray} Let
$\cW$
be a locally finite collection of subsets of
finite diameter in the ray
$[0,\infty)$.
Then there is a 1-Lipschitz
homeomorphism
$f:[0,\infty)\to [0,\infty)$
such that every set
$f(W)$, $W\in\cW$
has diameter at most 1.
\end{Lem}

\begin{proof} We let
$r_0=[0,\infty)$
and
$V_1$
be the set of all
$W\in\cW_0=\cW$
which intersect the initial segment of length 1 of
$r_0$.
Since the collection
$\cW$
is locally finite,
$V_1$
is finite. Then
$t_1=\sup\set{t\in r_0}{$t\in W\in V_1$}$
is finite, and we define
$\phi_1:r_0\to r_0$
by
$\phi_1(t)=t$
for
$t\in[0,1]$,
and for
$t\ge 1$
by
$\phi_1(t)=t$
if
$t_1\le 2$,
$\phi_1(t)=\frac{1}{t_1-1}(t-1)+1$
otherwise. Then
$\phi_1$
is a 1-Lipschitz homeomorphism, and the diameter
$\diam(\phi_1(W))\le 2$
for every
$W\in V_1$.

We put
$f_1=\phi_1$
and note that
$\cW_1=f_1(\cW)\cap r_1$
is a locally finite collection of subsets in the ray
$r_1=[1,\infty)$.
Now, we apply the same procedure to the ray
$r_1$
and the collection
$\cW_1$
extending the resulting 1-Lipschitz homeomorphism
$\phi_2:r_1\to r_1$
to
$r_0$
by the identity and putting
$f_2=\phi_2\circ f_1$.

Repeating we obtain a stabilizing sequence of 1-Lipschitz homeomorphisms
$f_n=\phi_n\circ f_{n-1}:r_0\to r_0$,
$f_n=f_{n-1}$
on
$[0,n]$,
and for which
$\diam(f_n(W))\le 2$
for all
$W\in\cW$
intersecting
$[0,n]$.
Composing the limit homeomorphism
$\lim_{n\to\infty}f_n$
with the homothety
$\la(t)=t/2$,
we obtain a required homeomorphism
$f$.
\end{proof}

\subsection{Cone over the boundary at infinity}

The following estimate from below is the major step in the
proof of Theorem~\ref{Thm:main}.

\begin{Pro}\label{Pro:hypdimcone} Let
$X\sub\hyp^n$, $n\ge 2$,
be the convex hull of a compact infinite
$Z\sub\di\hyp^n$.
Then
$\hypdim X\ge\dim Z+1$.
\end{Pro}

For the proof we need some preparations. We fix a base point
$o\in X$
and define the cone over
$Z$, $\cone(Z)\sub X$,
as the union of all geodesic rays emanating from
$o$
towards
$Z$.
Note that
$\cone(Z)$
with the metric induced from
$\hyp^n$
in general is neither CAT($-1$) nor even geodesic
space as for example in the case
$\dim Z=0$.
On the other hand,
$\cone(Z)$
is cobounded in
$X$
(see \cite[Proposition~10.1(2)]{BoS}), that is
$\dist(x,\cone(Z))\le\si_0$
for some
$\si_0>0$
and every
$x\in X$.
In what follows, we also use the angle metric in
$\di\hyp^n$
based at
$o$.

Next, we consider the annulus
$\an(Z)\sub\cone(Z)$,
which consists of all
$x\in\cone(Z)$
with
$1\le|x-o|\le 2$.
Clearly,
$\an(Z)$
is homeomorphic to
$Z\times I$, $I=[0,1]$.
According to a well known result from the dimension
theory (see \cite{Al}), the topological dimension
\[\dim\an(Z)=\dim Z+1 \tag{$\ast$}.\]
We also need the following simple fact from the dimension
theory.

\begin{Lem}\label{Lem:remove} Let
$Z$
be an infinite compact metric space. Then for every finite
$A\sub Z$
and all sufficiently small
$\ep>0$
we have
$\dim(Z\sm D_\ep(A))=\dim Z$,
where
$D_\ep(A)$
is the open
$\ep$-neighborhood
of
$A$.
\end{Lem}

\begin{proof} Since
$Z$
is infinite, one can assume that
$n=\dim Z>0$.
Let
$Y\sub Z$
be the subset which consists of all points
$z\in Z$
at which
$Z$
has dimension
$n$.
It is well known (see \cite[Ch.IV.5]{HW}) that
$\dim Y=n$,
thus
$Y$
cannot be finite, and
$Y\not\sub D_\ep(A)$
for all sufficiently small
$\ep>0$.
The claim follows.
\end{proof}

\begin{proof}[Proof of Proposition~\ref{Pro:hypdimcone}]
The space
$\hyp^n$
and hence its subspace
$X$
certainly have a bounded geometry. Let
$\rho_{\hyp^n}>0$
and
$M_{\hyp^n}:(0,\infty)\to(0,\infty)$
be the corresponding bounding parameters. Assume that
$\hypdim X<\dim Z+1$.
Then, moreover,
$\hypdim'X<\dim Z+1$
(see sect.~\ref{subsect:dmult}). Thus for
$d\ge\rho_{\hyp^n}$
there are a function
$N\in\cN$
and a uniformly
$N$-bounded
covering
$\cU$
of
$X$
with
$d$-multiplicity
$\le\dim Z+1$.
We fix
$\si\ge 10d$
and note that taking
$\si$-connected
components of
$\cU$
changes neither
$N$-boundedness
nor
$d$-multiplicity.
Thus we can
assume W.L.G. that the elements of
$\cU$
are
$\si$-connected.
By Lemma~\ref{Lem:boundend} there are only
finitely many elements
$U\in\cU$
with infinite diameter, and the boundary at
infinity each of them is finite.

We let
$\cV$
be the set of
$\si$-connected
components of the induced covering
$\cU\cap\cone(Z)$
of
$\cone(Z)$.
Then
$\cV$
is uniformly
$N$-bounded
and its
$d$-multiplicity
is
$\le\dim Z+1$.
The covering
$\cV$
of
$\cone(Z)$
can be represented as the disjoint union
$\cV=\cV_0\cup\cV_\infty$,
where
$\cV_0$
consists of all
$V\in\cV$
with finite diameter, and respectively
$\cV_\infty$
consists of all
$V\in\cV$
with infinite diameter.

There are natural polar coordinates
$x=(z,t)$, $z\in Z$, $t\ge 0$,
in
$\cone(Z)$.
Every 1-Lipschitz homeomorphism
$f:[0,\infty)\to[0,\infty)$
induces a 1-Lipschitz homeomorphism
$F:\cone(Z)\to\cone(Z)$
by
$F(z,t)=(z,f(t))$,
which does not change the visual diameter
$\diam(\sh(A))$
of any
$A\sub\cone(Z)$.

Given
$V\in\cV_0$
let
$W=W_V\sub[0,\infty)$
be the ray projection of
$V$,
i.e.,
$$W=\set{t\ge 0}{$(z,t)\in V\ \text{for some}\
    z\in Z$}.$$
Then
$\diam W<\infty$
and the collection
$\cW=\set{W_V}{$V\in\cV_0$}$
is locally finite in
$[0,\infty)$
since
$Z$
is compact and
$\cV_0$
is locally finite. By radial contraction Lemma~\ref{Lem:coveray}
there is 1-Lipschitz homeomorphism
$f:[0,\infty)\to[0,\infty)$
such that every set
$f(W)\sub[0,\infty)$, $W\in\cW$,
has diameter at most 1. For the induced homeomorphism
$F:\cone(Z)\to\cone(Z)$
we denote by
$\cV^1$
the covering
$F(\cV)=\set{F(V)}{$V\in\cV$}$
of
$\cone(Z)$.
Then
$\cV^1=\cV_0^1\cup\cV_{\infty}^1$
for
$\cV_0^1=F(\cV_0)$, $\cV_{\infty}^1=F(\cV_\infty)$,
and every
$V\in\cV_0^1$
lies in some annulus of width 1 centered at
$o$.

Consider the sequence of contracting homeomorphisms
$F_k:\cone(Z)\to\cone(Z)$
given by
$F_k(z,t)=(z,\frac{1}{k}t)$, $(z,t)\in\cone(Z)$, $k\in\N$,
and the corresponding sequence of coverings
$\cV^k=\cV_0^k\cup\cV_{\infty}^k$,
where
$\cV^k=F_k(\cV^1)$.
Using Lemma~\ref{Lem:ubgshadow} and the fact that
$F_k$
does not change the visual diameter, one easily obtains
that the diameter of the elements from
$\an(Z)\cap\cV_0^k$
vanishes as
$k\to\infty$.

On the other hand, the angle measure of
$U\sm B_c(o)$, $U\in\cU$,
is exponentially small in
$c$
by Lemma~\ref{Lem:ubgshadow}. It follows that for every
$\ep>0$
the shadow of
$U\sm B_c(o)$
is contained in the
$\ep$-neighborhood
of the finite set
$\di U$,
if
$c$
is chosen sufficiently large,
$\sh(U\sm B_c(o))\sub D_\ep(\di U)\sub\di\hyp^n$.
Thus for each
$V\in\cV_\infty$
the sequence
$F_k\circ F(V)\cap\an(Z)$
converges to a finite union of segments
$(z,[1,2])\sub\an(Z)$,
$z\in Z_V$,
as
$k\to\infty$,
where
$Z_V\sub\di U$
for
$U\in\cU$, $U\supset V$.
Since there are only finitely many elements
$U\in\cU$
with infinite diameter, and every
$V\in\cV_\infty$
is contained in one of them, there is a finite
$A\sub Z$
such that
$Z_V\sub A$
for every
$V\in\cV_\infty$.
Moreover, for every
$\ep>0$
there is
$k_\ep\in\N$
such that the compact set
$\an_\ep(Z)=\an(Z\sm D_\ep(A))$
misses any
$V\in\cV_{\infty}^k$,
and thus it is covered by elements of
$\cV_0^k$, $k\ge k_\ep$.
The multiplicity of the covering
$\cV_0^k$
of
$\an_\ep(Z)$
is
$\le\dim Z+1$,
and the diameter of its elements vanishes as
$k\to\infty$.
Thus
$\dim\an_\ep(Z)\le\dim Z$.
However, by ($\ast$) and Lemma~\ref{Lem:remove},
$\dim\an_\ep(Z)=\dim Z+1$
for all sufficiently small
$\ep>0$.
This is a contradiction.
\end{proof}

\subsection{Proofs of Theorem~\ref{Thm:main} and Corollary~\ref{Cor:nonemb}}

\begin{proof}[Proof of Theorem~\ref{Thm:main}] By \cite[Theorem~1.1]{BoS},
the space
$X$
is roughly similar to a convex subset of
$\hyp^n$
for some integer
$n$.
Actually, it is proved there that
$\di X$
is homeomorphic to a compact
$Z\sub\hyp^n$
such that
$X$
is roughly similar to the convex hull of
$Z$
in
$\hyp^n$.
Thus we can assume that
$X$
is the convex hull of some compact, infinite
$Z\sub\hyp^n$.
Now, Proposition~\ref{Pro:hypdimcone} completes the proof.
\end{proof}

\noindent
{\em Proof of Corollary~\ref{Cor:nonemb}.}
We actually prove a more general result.
\begin{Thm}\label{Thm:nonemb} Assume that there is
a quasi-isometric embedding
$$f:X\to T_1\times\dots\times T_k\times Y,$$
of a Gromov
hyperbolic space
$X$,
satisfying the conditions of Theorem~\ref{Thm:main}, into
the
$k$-fold
product of metric trees
$T_1,\dots,T_k$
stabilized by an UBG-factor
$Y$.
Then
$k\ge\dim\di X+1$.
\end{Thm}

\begin{proof} First, we note that
$\hypdim Y=0$.
Next, the asymptotic dimension of every metric tree
$T$
is at most 1,
$\asdim T\le 1$,
see \cite[Proposition~4]{DJ}.
Therefore,
$\hypdim T\le 1$,
and by the Product Theorem the hyperbolic dimension of
the target space is at most
$k$.
On the other hand,
$\hypdim X\ge\dim\di X+1$
by Theorem~\ref{Thm:main}. Now, the required estimate follows
from mononicity of
$\hypdim$,
see Theorem~\ref{Thm:monhypdim}.
\end{proof}


\bigskip
\begin{tabbing}

Sergei Buyalo,\hskip11em\relax \= Viktor Schroeder,\\

St. Petersburg Dept. of Steklov \>
Institut f\"ur Mathematik, Universit\"at \\

Math. Institute RAS, Fontanka 27, \>
Z\"urich, Winterthurer Strasse 190, \\

191023 St. Petersburg, Russia\>  CH-8057 Z\"urich, Switzerland\\

{\tt sbuyalo@pdmi.ras.ru}\> {\tt vschroed@math.unizh.ch}\\

\end{tabbing}

\end{document}